\documentclass[journal]{IEEEtran}
\usepackage[cmex10,intlimits]{amsmath}
\usepackage{amsfonts}
\usepackage{amssymb}
\usepackage{graphicx}
\usepackage{booktabs}
\usepackage{lipsum}
\usepackage{xcolor}
\usepackage{tikz}
\usetikzlibrary{hobby}
\usepackage{cite}
\usepackage[english]{babel}

\usepackage[most]{tcolorbox}

\newtcolorbox[auto counter]{CMAbox}[2][]{before upper={\parindent15pt\noindent},floatplacement=h!,float,center title,%
colback=PairedA!10,colframe=PairedB,fonttitle=\bfseries,fontupper=\small,title=Box~\thetcbcounter. #2,#1}

\usepackage{capekCommands}
\usepackage{capekColors}

\graphicspath{figures/}

\begin{document}

\bstctlcite{IEEEexample:BSTcontrol}
\title{Computational Aspects of Characteristic Mode Decomposition -- An Overview}
\author{Miloslav~Capek,~\IEEEmembership{Senior Member,~IEEE}
        and~Kurt~Schab,~\IEEEmembership{Member,~IEEE}
\thanks{Manuscript received \today; revised \today. This work was supported by the Czech Science Foundation under project~\mbox{No.~21-19025M}.}
\thanks{M. Capek is with the Department of Electromagnetic Field, Faculty of
Electrical Engineering, Czech Technical University in Prague, 166 27 Prague,
Czech Republic (e-mail: miloslav.capek@fel.cvut.cz).}
\thanks{K. Schab is with the Department of Electrical and Computer Engineering, Santa Clara University, Santa Clara, CA, USA (e-mail: kschab@scu.edu).}
\thanks{Manuscript received April 19, 2005; revised September 17, 2014.}
}

\maketitle

\begin{abstract}
Nearly all practical applications of the theory of characteristic modes (CMs) involve the use of computational tools.  Here in Paper 2 of this Series on CMs, we review the general transformations that move CMs from a continuous theoretical framework to a discrete representation compatible with numerical methods.  We also review key concepts encountered across a variety of numerical CM implementations. These include modal tracking, dynamic range, code validation, and techniques related to electrically large problems.
\end{abstract}


\IEEEpeerreviewmaketitle
\section{Introduction}
\label{sec:intro}

\IEEEPARstart{W}{hile} the theory of characteristic modes (CMs) can be applied rigorously in certain theoretical cases, it has long been associated with computational methods allowing for the analysis of arbitrarily shaped systems.  The method of moments (MoM) \cite{Harrington_FieldComputationByMoM} in particular is closely tied to the analysis of CMs due to its direct transformation of the continuous impedance operator into a finite-dimensional impedance matrix \cite{Harrington_1971b}.  Simple MoM implementations are commonly taught in graduate-level coursework on computational electromagnetics \cite{Makarov_AntennaAndEMModelingWithMatlab,jin2011theory, Harrington_FieldComputationByMoM}, and many commercial software packages allow for the generation and exporting of MoM-related data \cite{atom,feko}.  Contemporary computing power allows for the solution of moderately sized MoM problems (on the order of a few thousand unknowns) without the necessity of employing sophisticated algorithms and expensive high-performance hardware.  This availability of MoM codes and their close connection to modal analyses has enabled many researchers to explore different aspects of CM theory without excessive overhead spent in the development of specialized computational tools.  As such, many studies in the literature implicitly assume MoM numerical tools as the basis for the theory of CMs, rather than the continuous operator theory underpinning its original formulations~\cite{Garbacz_TCMdissertation}.  While this is not entirely accurate, it does reflect the prevalence of numerical tools in the study of CMs and motivates the present second part of this review series on computational aspects of characteristic mode analysis (CMA).

The intent of this review is to highlight key challenges and techniques in the numerical implementation of CMA, with focus placed on methods preserving the favorable properties associated with purely analytical formulations. This part follows general overview (Paper 1 of this series) and precedes a review on canonical antenna design with the help of CMs~(Paper 3 of this series).

\section{Characteristic Modes}

CM decomposition is motivated by the desire for a modal basis with properties convenient for electromagnetic scattering analyses.  Most often, these types of decomposition are based on the diagonalization of specific operators arising from frequency domain integral equations. Beginning in the 1940's, decompositions of fields scattered by obstacles were independently introduced using the scattering matrix~\cite{1948_Montgomery_Principles_of_Microwave_Circuits}, the transition matrix \cite{Garbacz_TCMdissertation}, and the impedance matrix \cite{Harrington_1971a}. An important property of the CM decomposition is that the external (radiating) problem is treated, making it possible to orthogonalize modal radiation patterns.

Contemporary work on CMs is concentrated on formulations based on impedance operators relating current distributions on a specified object to the fields they scatter, see Box~1.  In this case, CMs are commonly defined as those current distributions diagonalizing both the impedance and radiation operators. Modifications of the CM decomposition exist, see Section~\ref{sec:reltech}, though these techniques often come at the cost of certain properties which make the impedance / radiation formulation favorable~\cite{Harrington_1972b}.

CMs represent a basis of continuous current distributions and are defined independently of any computational method.  However, outside of select canonical problems, \eg{}, \cite{SarkarMokoleSalazarPalma_AnExposeOnInternalResonancesCM,Bernabeu_Jimenez_2017a,Huang_StudyontheRelationshipsbetweenEigenmodesNaturalModesandCharacteristicModes}, solutions to the CM eigenvalue problem can only be obtained through the use of numerical tools. 

\section{Matrix Representation of CMs}
\label{sec:arbshapes}

Computational implementation of the theory of CMs relies on the conversion of fields and operators defined over continuous domains into vectors and matrices corresponding to a basis with discrete elements~\cite{Harrington_1971b}.  By approximating currents on arbitrarily shaped conducting surfaces through an appropriate basis (\eg{}, RWG basis functions~\cite{RaoWiltonGlisson_ElectromagneticScatteringBySurfacesOfArbitraryShape}, see top insets in Box~\ref{A:Box1}) the electric field integral equation (EFIE) may be transformed into a system of equations via application of MoM~\cite{Harrington_FieldComputationByMoM}. This conversion, in turn, leads to a discrete representation of the CM eigenvalue problem, see Box~1 and bottom insets therein. 

When the impedance matrix is transpose symmetric, \eg{}, when Galerkin method~\cite{Harrington_FieldComputationByMoM} is used within the MoM~to model reciprocal systems, CMs are equiphase and diagonalize the impedance operator, leading to orthogonality in both the radiation and reactance operator. Normalization of CMs is arbitrary, though normalizations forcing each modal current to radiate equal power are by far the most common, see Box~2.  

In some cases, computation of the impedance matrix may be accelerated via non-Galerkin procedures or asymmetric integrations \cite{Makarov_AntennaAndEMModelingWithMatlab},~leading to asymmetric impedance matrices. In such cases care must be taken to either be aware of reduced orthogonality or to artificially \Quot{symmetrize} appropriate matrix operators.  Orthogonality in the matrix $\M{R}$, by Poynting's theorem, also leads to orthogonality in modal far-fields when integrated over the entire far-field sphere, see Box~2.  Note that this orthogonality property is not, in general, present for lossy structures unless modifications to the generating eigenvalue problem are altered, necessarily impacting other orthogonality properties in the process~\cite{Harrington_1972b,Chang_1977a, Yla_Oijala_2019a, Yla_Oijala_2019b}.  Additionally, modal currents are not, in general, orthogonal in an unweighted inner product, \ie{}, $\M{I}^\T{H}_m\M{I}_n \neq c_{mn}\delta_{mn}$, and their orthogonality is reserved only for so-called separable coordinate systems~\cite{morse1953methods} for which analytical solutions to the problem~\eqref{eq:MR3} may typically be found~\cite{Stratton_ElectromagneticTheory, Garbacz_TCMdissertation, SarkarMokoleSalazarPalma_AnExposeOnInternalResonancesCM, Amendola1997, Capek_2017b}.
The orthogonality relations in \eqref{eq:NM1} and \eqref{eq:NM2} imply modal orthogonality in reactive energy, \ie{}, the difference between stored magnetic and electric energies.  While the reactance operator~$\M{X}$ can be decomposed into electric and magnetic energy operators \cite{Vandenbosch_ReactiveEnergiesImpedanceAndQFactorOfRadiatingStructures, GustafssonTayliEhrenborgEtAl_AntennaCurrentOptimizationUsingMatlabAndCVX}, orthogonality of modes does not, in general, extend to these individual operators or their summation which represents stored electromagnetic energy.  

Due to their relative excess of magnetic (electric) energy, modes with positive (negative) eigenvalues are referred to as inductive (capacitive), \cf{} \eqref{eq:NM2} in Box~2.  A balance of modal electric and magnetic energies leads external resonance \footnote{Here \emph{external resonance} denotes the condition $\lambda_n = 0$ whereas \emph{internal resonances} are typically associated with non-radiating currents.} at $\lambda = 0$.  The spectrum of eigenvalues spans all positive and negative real numbers, with accumulation at $\pm\infty$.  As such, eigenvalues are frequently reported in terms of log magnitude or characteristic angle~\cite{Newman_SmallAntennaLocationSynthesisUsingCharacteristicModes} $\alpha_n = \pi - \arctan \lambda_n$, the latter being restricted to the interval $[\pi/2,3\pi/2]$.  Similarly the quantity \mbox{$0\leq|1+\T{j}\lambda_n|^{-1}\leq 1$} is frequently used as a metric of modal significance due to its predictable range and appearance in many CM expansion formulae.

\begin{CMAbox}[label={A:Box1}]{{EFIE Formulation and Its Matrix Representation}}

{\centering
  \includegraphics[width=\textwidth]{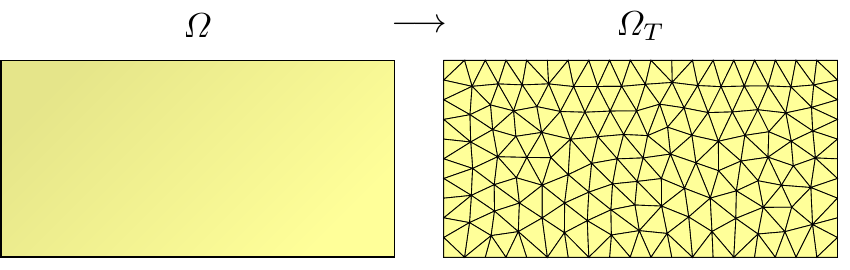}\vspace{0.25cm}\par}

\noindent
The electric field integral equation (EFIE) relates a scattered field~$\V{E}_\T{s} (\V{r})$ and its sources~$\V{J}(\V{r}')$ as~\cite{Harrington_TimeHarmonicElmagField}
\begin{equation}
\V{E}_\T{s} (\V{r}) = - \J \omega \MUE \int\limits_\srcRegion \M{G} (\V{r},\V{r}') \cdot \V{J}(\V{r}') \D{\V{r}'}
\label{eq:MR1}
\end{equation}
with $\M{G} (\V{r},\V{r}')$ being the dyadic Green's function.
The above relation can be use to define an impedance operator~\mbox{$\OP{Z}(\V{J}) = \OP{R}(\V{J}) + \J \OP{X}(\V{J}) = \UV{n} \times \UV{n} \times \V{E}_\T{s} (\V{J})$} from which a CM eigenvalue problem is constructed as~\cite{Harrington_1971a}
\begin{equation}
\OP{Z} ( \V{J}_n) = \left( 1 + \J \lambda_n \right) \OP{R} ( \V{J}_n).
\label{eq:MR3}
\end{equation}
The generalized eigenvalue problem~\eqref{eq:MR3} demands a solution to the integral equation, which can be analytically found in exceptional cases only. For this reason, the continuous quantities are represented in a suitable basis~$\left\{\basisFcn_m(\V{r}')\right\}$,~\cite{PetersonRayMittra_ComputationalMethodsForElectromagnetics}
\begin{equation}
\V{J} (\V{r}') \approx \sum_{n=1}^N I_n \basisFcn_n (\V{r}'), \quad \V{r}' \in \srcRegion_T,
\label{eq:MR4}
\end{equation}
which yields a matrix representation of the impedance operator $\M{Z} = \M{R} + \J \M{X} \equiv [Z_{pq}]$ with elements given by
\begin{equation}
Z_{pq} = \J \omega \MUE \int\limits_{\srcRegion_T} \int\limits_{\srcRegion_T} \basisFcn_p(\V{r}) \cdot \M{G} (\V{r},\V{r}') \cdot \basisFcn_q(\V{r}') \D{\V{r'}} \D{\V{r}}.
\label{eq:MR5}
\end{equation}
From this representation, a matrix form of the CM decomposition~\cite{Harrington_1971b}
\begin{equation}
\M{Z} \Ivec_n = \left( 1 + \J \lambda_n \right) \M{R} \Ivec_n\quad\leftrightarrow\quad\M{X}\Ivec_n = \lambda_n\M{R}\Ivec_n,
\label{eq:MR6}
\end{equation}
may be constructed. It should be noted that the choice of basis~\eqref{eq:MR4} is arbitrary, \ie{}, it can be composed of entire-domain (\eg{}, spherical waves) or piece-wise defined (\eg{}, RWG) basis functions. While the transition from analytical~\eqref{eq:MR3} to numerical~\eqref{eq:MR6} representation of the CM problem significantly broadens the scope of its practical use, it also restricts the inherent precision and may lead to numerical issues, discussed further in Section~\ref{sec:dr}.

{\vspace{0.25cm}\centering
  \includegraphics[width=\textwidth]{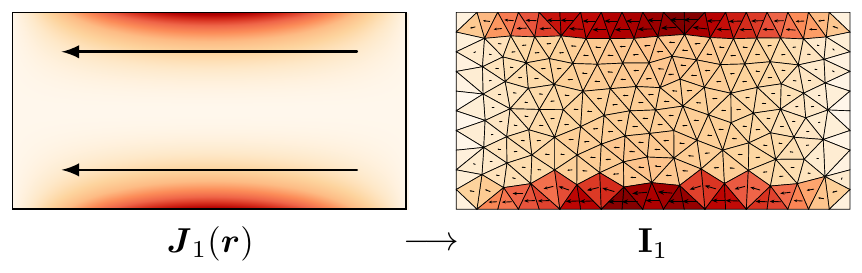}\par}

\end{CMAbox}

\begin{CMAbox}[label={A:Box2}]{{Current Normalization}}

\noindent
Transpose symmetry of the impedance matrix~$\M{Z}$ leads to CM current distributions which are orthogonal in both the radiation and reactance operators. The orthogonality in the radiation operator implies orthogonal far fields (see the inset below). The choice of normalization, however, is arbitrary and provides an opportunity to re-scale eigenvectors~$\V{J}_n(\V{r}')$ and~$\Ivec_n$ at will. The most commonly adopted normalization ensures real-valued current vectors radiating equal amounts of power, \ie{},
\begin{multline}
P_{\T{r},n} = \dfrac{1}{2 Z_0} \int_{4\pi} \V{F}_m^\ast (\UV{r}) \cdot \V{F}_n (\UV{r}) \D{\varOmega} \\ = \dfrac{1}{2} \langle \V{J}_m, \OP{R} (\V{J}_n) \rangle \approx \dfrac{1}{2} \Ivec_m^\herm \M{R} \Ivec_n = \delta_{mn}.
\label{eq:NM1}
\end{multline}
This choice of normalization also equates modal reactive power with the CM as
\begin{multline} 2\omega\left(W_{\T{m},n}-W_{\T{e},n}\right) = \dfrac{1}{2} \langle \V{J}_m, \OP{X} (\V{J}_n) \rangle\\ \approx \dfrac{1}{2} \Ivec_m^\herm \M{X} \Ivec_n = \lambda_n \delta_{mn},
\label{eq:NM2}
\end{multline}
where~$W_{\T{m},n}$ and~$W_{\T{e},n}$ are the modal magnetic and electric stored energies~\cite{2018_Schab_Wsto}.

As an example, the radiation pattern cuts of the first four characteristic modes of a perfectly electric conducting plate from Box~\ref{A:Box1} are depicted for~$ka = 1$. When~\eqref{eq:NM1} is applied, their size is changed dramatically. 

{\vspace{0.25cm}\centering
  \includegraphics[width=\textwidth]{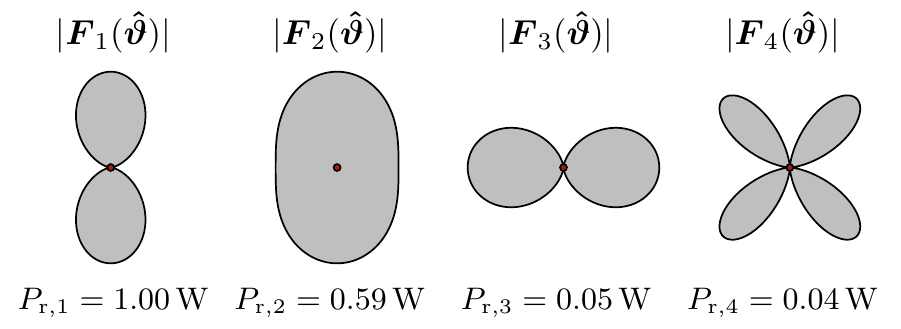}\par}

\end{CMAbox}

\section{Accuracy, Dynamic Range, and Convergence}
\label{sec:dr}

Solutions to the discrete eigenvalue problem in \eqref{eq:MR6} may be found using a variety of numerical methods.  The choice of numerical method, however, along with the properties of the underlying impedance matrix data, may lead to noticeable differences in modal solutions.  In this section, we discuss key numerical aspects of obtaining CM data.

\subsection{Solving Generalized Eigenvalue Problems}

When it comes to the numerical evaluation of~\eqref{eq:MR6}, two algorithms are typically employed: generalized Schur decomposition (GSD) and the implicitly restarted Arnoldi method (IRAM)~\cite{golub2013matrix}. For example, in MATLAB~\cite{matlab}, their LAPACK and ARPACK implementations are accessible via \texttt{eig()} and \texttt{eigs()} functions, respectively, and similarly for other high-level programming languages. The most notable difference is their algorithmic complexity. The GSD method achieves higher numerical stability and produces a full set of eigenmodes, however, the cost is $\OP{O}(N^3)$, where~$N$ is the number of degrees of freedom, \cf{}, \eqref{eq:MR3}. The IRAM method is capable of efficiently finding only the first~$M$ modes with complexity~$\OP{O}(M N^2)$. When the eigenvalue problem is solved, the CMs are typically then normalized, see Box~2, and tracking is performed, see Box~3.

\subsection{Modal Dynamic Range}
For a system with~$N$ degrees of freedom, the generalized eigenvalue problem~\eqref{eq:MR6} defines a set of~$N$ CMs.  However, issues related to numerical precision (specifically in the representation of weakly radiating currents with complex spatial dependence) typically limit meaningful results to those \mbox{$M\ll N$} modes with eigenvalue magnitudes below an empirically observed threshold.  Conveniently, \Quot{low-order} modes with eigenvalues at or near resonance are those most commonly studied \cite{Antonino_Daviu_2006a,Adams_2011a,Ghosal_2020a}.  The numerical threshold determining the number of available modes depends on many factors, particularly singularity treatment and numerical precision \cite{Capek_2017b}, but also the operational frequency and object complexity.

\begin{figure}
\centering
\includegraphics[]{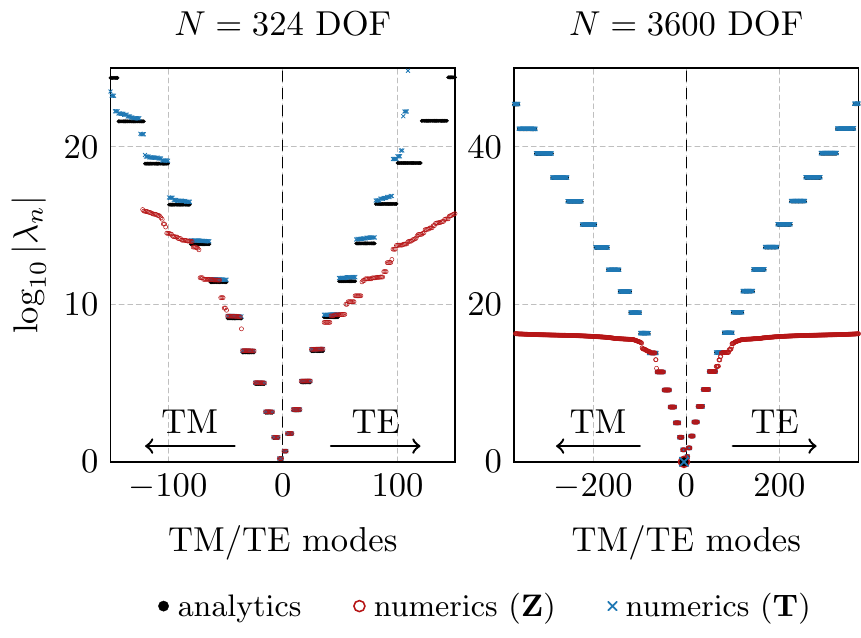}
\caption{CM decomposition of a spherical shell of radius $a$ discretized into $216$ (left) and $2400$ (right) triangles, leading to $324$ and $3600$ basis functions, respectively. Electrical size of the spherical shell is $ka=1$. Analytical results~\cite{Capek_2017b} (black dots) are compared with two numerical methods based on the decomposition of the impedance matrix~$\M{Z}$ (red circles)~\eqref{eq:MR6} and transition matrix~$\M{T}$ (blue $\times$)~\cite{Garbacz_TCMdissertation}. See~\cite{Capek_2017b} for comparison with commercial software.}
\label{fig:sphere}
\end{figure}

An important issue significantly affecting the numerical dynamics of the CM generalized eigenvalue problem (CMGEP) is the null-space of the $\M{R}$ matrix~\cite{Harrington_1971b}. This null-space can be extracted, which greatly extends the dynamic range of available modes and reduces computational time~\cite{Tayli_2018a}. An ultimate solution in this direction is to evaluate the transition matrix~$\M{T}$ from the impedance matrix~$\M{Z}$ as described in~\cite{Losenicky_etal_MoMandThybrid} and perform decomposition in a basis of outgoing spherical waves. This leads to the improved performance shown in Fig.~\ref{fig:sphere}. Another problem is low-frequency instability which is addressed in~\cite{Dai_2017a}.

Figure~\ref{fig:dipole-modes} demonstrates the limited number of CMs by plotting modal current distributions (top) and eigenvalue magnitude (bottom left) for a straight dipole antenna. There it is clear that for modes with eigenvalue magnitudes less than $\sim10^{16}$, modal current distributions are smooth functions resembling a Fourier basis in one dimension, whereas modes with eigenvalues saturated above this level have current distributions that are effectively numerical noise.

\subsection{Convergence in Modal Superposition}

\begin{figure}
\centering
\includegraphics[width=3.15in]{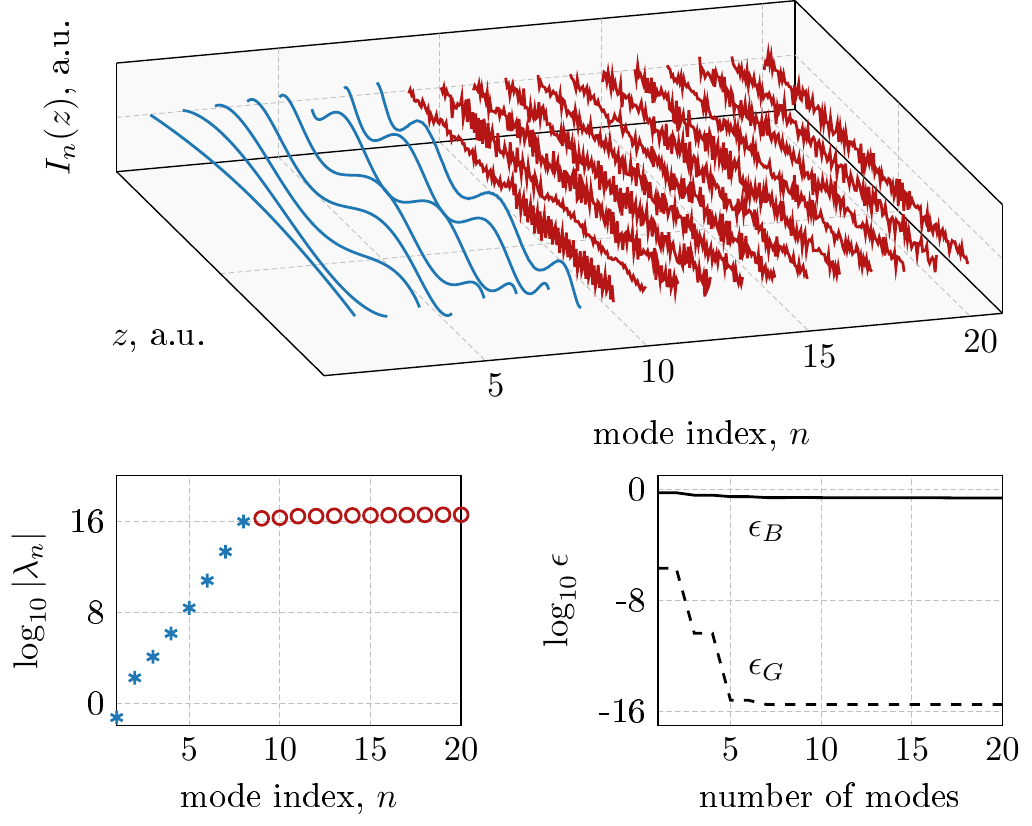}
\caption{CM properties of a wire dipole.  Modes are ordered by increasing eigenvalue magnitude.  Modal current distributions (top) and eigenvalue magnitudes (bottom left) are shown with blue and red markers representing modes below and above the numerical threshold, respectively.  Normalized residual input conductance and susceptance (bottom right) indicates poor convergence.}
\label{fig:dipole-modes}
\end{figure}

Diagonalization of the impedance operator leads to a simple form of expansion coefficients for the current induced by an incident field represented by the vector~$\M{V}$,
\begin{equation}
\M{I} = \sum_n \alpha_n \M{I}_n,\quad \alpha_n = \frac{\M{I}^\T{H}_n\M{V}}{2 P_{\T{r},n}\left(1+\T{j}\lambda_n\right)},
\end{equation}
where~$\alpha_n$ is commonly referred to as a modal weighting coefficient. The explicit form of the modal weighting coefficients suggests that modal contributions to the total current distribution depend both on the alignment (inner product) of modal currents with the excitation ($\M{I}_n^\T{H}\M{V}$) as well as the reactive nature of each mode (via $|1+\T{j}\lambda_n|^{-1}$). These two terms are commonly called the modal excitation coefficient and modal significance, respectively. Note that highly reactive currents with~$|\lambda_n|\gg 1$ may make a non-trivial contribution to the above sum due to large current magnitudes induced by the chosen normalization enforcing all currents to radiate unit power.  This is particularly true for the common, but poorly converging, case of antennas driven by delta gap sources~\cite{yee1973self}, as illustrated by the example data in the bottom right panel of Fig.~\ref{fig:dipole-modes} showing relative errors in driven input conductance ($\epsilon_G$) and susceptance ($\epsilon_B$) as a function of increasingly higher number of modes used to reconstruct the driven current distribution of a center fed dipole.  In synthesizing equivalent circuit models of driven antennas, excess reactance due to modes far away from resonance are frequently modeled using empirically determined reactive components, \eg{}, \cite{adams2013broadband}.

\section{Tracking}
\label{sec:tracking}

\begin{CMAbox}[label={A:Box3}]{{Crossings and Crossing Avoidances}}

{\vspace{0.25cm}\centering
  \includegraphics[width=\textwidth]{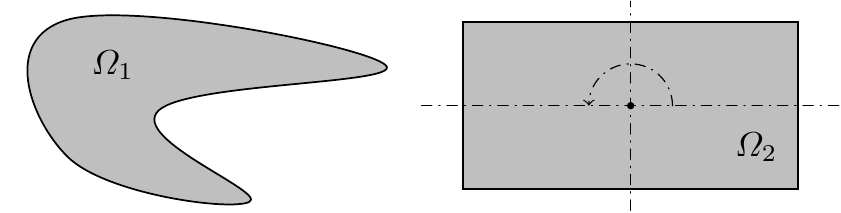}\par}
\vspace{0.2cm}

\noindent
The von Neumann-Wigner theorem~\cite{vonNeumannWigner_UberDasVerhaltenEigenvertwen} gives rise to several important results related to modal tracking and symmetry.  

In greatly simplified terms, this theorem states that there can be no eigentrace crossings so long as an object under study has no point-group symmetries (see the ``asymmetric'' shape~$\srcRegion_1$ above and the left inset below), \ie{}, the impedance matrix is in its irreducible form~\cite{Knorr_1973_TCM_symmetry}. On the contrary, when an object possesses point-group symmetries (\eg{}, mirror or rotational symmetry, see shape~$\srcRegion_2$ above), its impedance matrix is reducible into a block-diagonal matrix and eigentraces from different blocks (associated with irreducible representations) may cross each other, see the right inset below. The importance of this result extends far beyond CMs as it applies to any eigenvalue problem, \eg{} those encountered in quantum mechanics~\cite[Figs. 30, 31]{landau2013quantum}. See~\cite{Schab_2017a} and~\cite{Masek_2020a} for further details.

{\vspace{0.25cm}\centering
  \includegraphics[]{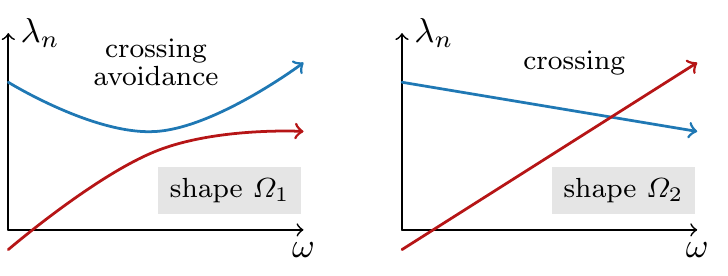}\par}

\end{CMAbox}

\begin{figure}
    \centering
    \includegraphics[width=3.5in]{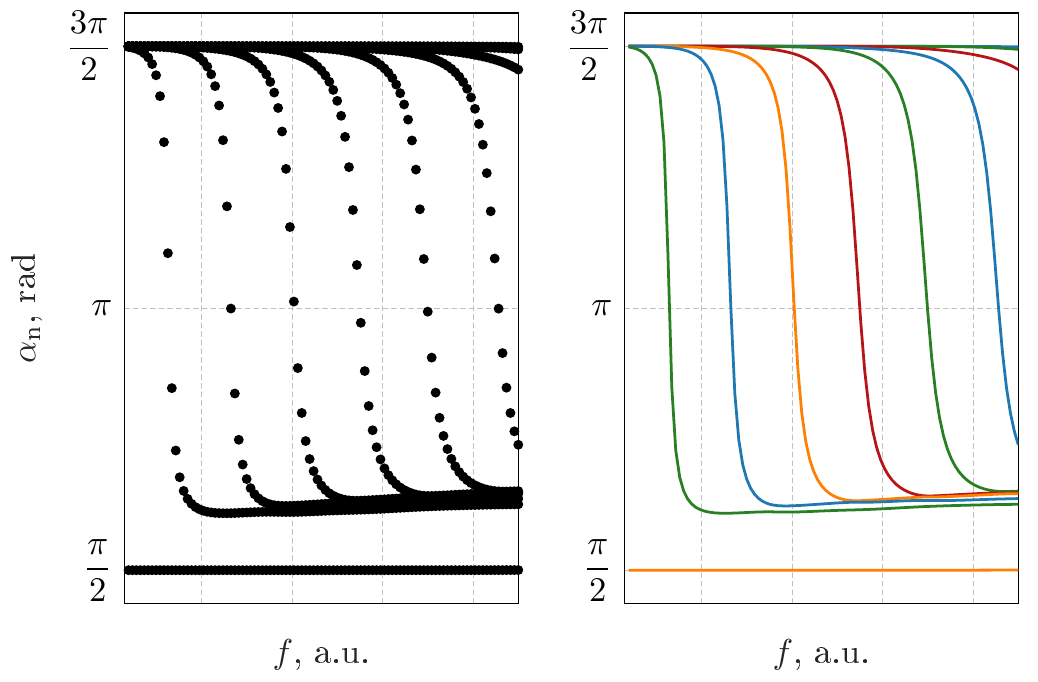}
    \caption{Untracked (dots, left) and tracked (lines, right) modal eigenvalues using analytic symmetry-based tracking rules for a structure with $C_1$ symmetry (no modal crossings allowed).}
    \label{fig:tracking}
\end{figure}

The operators involved in the CMGEP are naturally parameterized by frequency~$\omega$, leading to frequency dependent modal quantities.  The concept of associating modal quantities across frequencies (particularly over a set of discrete frequencies at which the CMGEP is solved) is known as \emph{modal tracking}.  Modal data represented by characteristic angles~$\alpha_n$ are shown schematically in Fig.~\ref{fig:tracking}, where untracked data at discrete data frequency points are tracked to create smooth eigentraces.

Based on hypotheses of various metrics of modal similarity, many heuristic tracking algorithms have been proposed, such as leveraging vector correlation \cite{Raines_2012a, Safin_2016a}, eigenvalue proximity \cite{Chen_2021a}, far-field correlation \cite{Miers_2015a}, and hybrid methods \cite{LudickJakobusVogel_AtrackingAlgorithmForTheEigenvectorsCalculatedWithCM}.  Complicating the process of modal tracking are degeneracies and crossing avoidances (associated with modal coupling) \cite{KingAJ_PhDThesis, Schab_2016a, Lin_2017a, Ghosal_2020b} near frequencies where multiple eigenvalues approach the same value. Though a system’s symmetry \cite{Knorr_1973_TCM_symmetry, Vescovo_1989a, Vescovo_1990a} may be employed to analytically distinguish degeneracies and crossing avoidances \cite{Schab_2017a} and enable smooth tracking of complex modal systems \cite{Masek_2020a}, the sensitivity of these phenomena to vanishingly small symmetry-destroying perturbations raises questions regarding the physical relevance of tracking modes across frequency.  Regardless of physical interpretation, one special case stands out among results related to symmetry-based tracking: modal eigenvalues of completely asymmetric ($C_1$ symmetry) structures are infinitely unlikely to cross \cite{vonNeumannWigner_UberDasVerhaltenEigenvertwen, Schab_2017a}.  Additionally, modal symmetry has factored into the design of MIMO systems with orthogonal radiation channels~\cite{Peitzmeier_2019a,Masek_2020a}.

Several misconceptions regarding modal tracking are prevalent in the literature.  First, it is not generally appropriate to assume \Quot{slowly changing} modal currents as functions of frequency.  As observed in other branches of physics involving coupled modes \cite{haus1991coupled}, modal eigenvectors may change rapidly in the vicinity of crossing avoidances \cite{Schab_2016a}.  Second, though useful for interpretation and visualization, the premise of modal tracking is primarily one of human convenience, rather than physical truth.  Symmetry-based tracking rules demonstrate that infinitesimal perturbations to a system lead to differing tracking results \cite{Masek_2020a}.

\section{Extended Techniques}
\label{sec:reltech}

So far in this review only the specific generalized eigenvalue problem in \eqref{eq:MR3}, formulated for the study of lossless scatterers in free space, has been considered.  However, many variations on this problem exist, allowing for CMs to be applied in a number of alternate settings.  

\subsection{Electrically Large Problems}
\label{sec:ellarge} 

\begin{figure}
\centering
\includegraphics[scale=0.84]{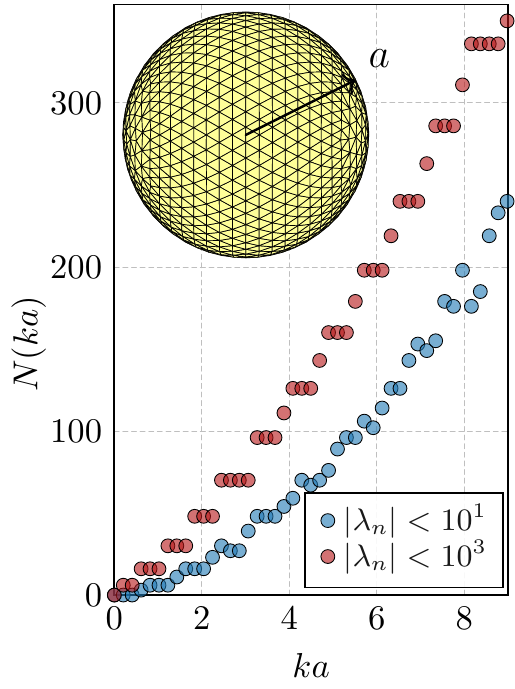}
\includegraphics[scale=0.84]{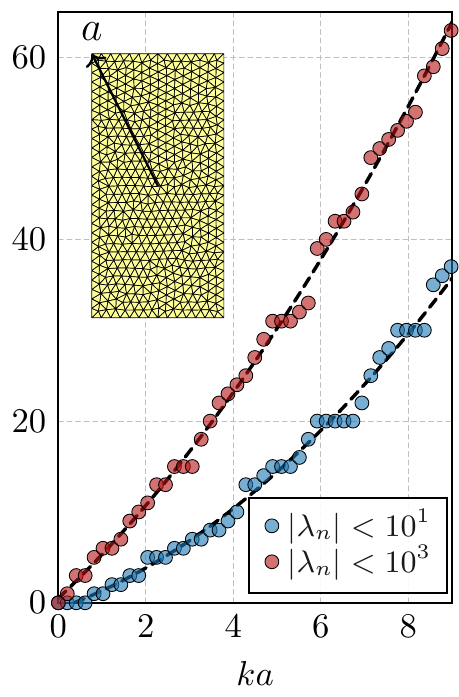}
\caption{Number of significant modes for a spherical shell (left) and a rectangular plate (right) depending on electrical size~$ka$. The fits for the rectangular plate are found for a second-degree polynomial and depicted by the dashed lines. Explicitly, the fits are $N(ka) \approx 0.26 (ka)^2 + 1.7 ka$ (blue markers) and $N(ka) \approx 0.28 (ka)^2 + 4.6 ka$ (red markers).}
\label{fig:numOfModes}
\end{figure}

As the electrical size of a structure increases, so does the number of degrees of freedom required to accurately model induced currents.  This increase, in turn, leads to rapidly growing computational cost in populating, storing, and manipulating an impedance matrix generated by EFIE-based MoM, limiting the practical feasiblity of carrying out CM analyses.  Fast multipole methods have been applied to remedy this issue, greatly reducing the computational cost involved in obtaining dominant CMs from large structures \cite{Dai_2016a,Wu_2018a,Zhang_2020b}. Nevertheless, it is important to point out that the number of significant modes increases rapidly with electrical size~$ka$, see Fig.~\ref{fig:numOfModes}.

\subsection{Dual and Substructure Problems}
\begin{figure}
    \centering
    \includegraphics[width=3.5in]{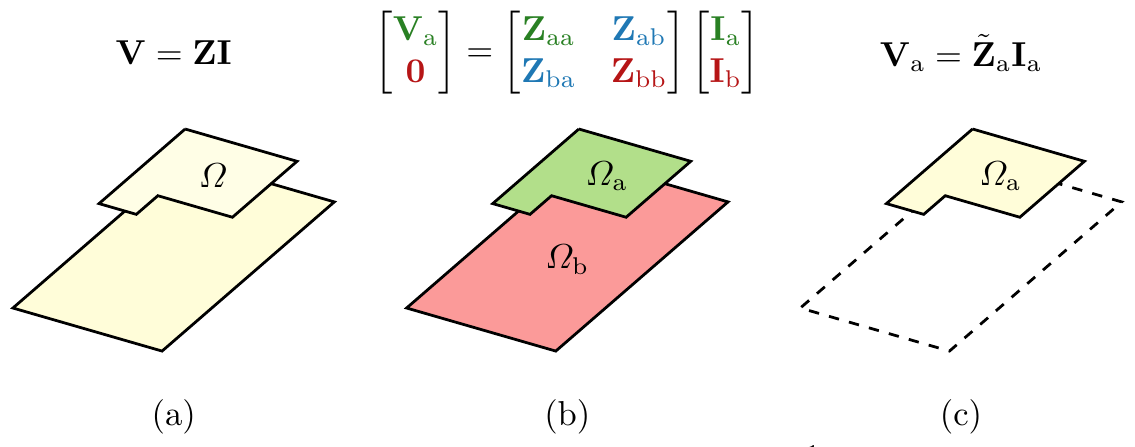}
    \caption{The problem of an antenna near a scattering object~$\srcRegion$ cast in three forms: (left) the entire structure represented as one impedance matrix system, (middle) the structure partitioned into antenna ($\srcRegion_\T{a}$) and scatterer ($\srcRegion_\T{b}$) regions, and (right) the antenna region alone represented by a compressed impedance matrix eliminating dependent scatterer currents.  Note that this example assumes no direct excitation on the scatterer, which is common in embedded and platform-based antenna designs.}
    \label{fig:substructure}
\end{figure}

Dual formulations involving Babinet’s principle may be employed to represent apertures and slot structures on infinite ground planes as magnetic currents, see \cite{Kabalan_2001a,Liang_2018b} and references therein.  Numerical Green’s functions may be employed to generalize the impedance operator to non-free-space settings, most commonly in the analysis of embedded antennas located on or near conducting platforms \cite{chalas2012efficient,Ethier_2012a,Ethier_2014a, Alroughani_2015a,Alakhras_2020a}.  When all objects in such a system are of finite extent, calculation of the numerical Green’s function is closely related to the partitioning of an impedance matrix representing the entire system, see Fig.~\ref{fig:substructure}, while problems involving infinite ground planes may be approached through the use of image theory \cite{Angiulli_2000a, Borchardt_2019a}.  

\subsection{Arrays and Network Modes}
The finite basis impedance matrix is equivalent to that encountered in the analysis of $N$-port networks, leading directly to the concept of network CMs \cite{HarringtonMautz_PatternSynthesisForLoadedNportScatterers, Obeidat_2010b, Jaafar_2017a}.  For array problems, the scattering problem may be reduced to an~$N$-port network or whole domain basis functions may be employed to approximate currents over each element to reduce computational overhead~\cite{Tzanidis_2012a, KingAJ_PhDThesis, Lonskyetal_CMofDipoleArrays,Guan_2019a}.  More involved approaches have been developed for similar analysis of finite arrays through the use of modal properties of single array elements \cite{Cheng_2019a}.  Numerical analysis of infinite uniform arrays has also been carried out through the use of problem-specific Green's functions, closely resembling the approaches taken in the analysis of objects above infinite ground planes \cite{Angiulli_2000a,maalik2016characteristic, Haykir_2019a}.  

\subsection{Related Modal Methods}
Variations on the CM eigenvalue problem involving weighted far-field integration (as opposed to the implicit uniform weighting present in the classical CMGEP) give rise to Inagaki modes \cite{InagakiGarbacz_EigenfunctionsOfCompositeHermitianOperatorsWithApplicationToDiscreteAndContunuousRadiatingSystems} and generalized CMs \cite{Liu_SomeRelationshipsBetweenCharacteristicModesAndInagakiModesForUseInScatteringAndRadiationProblems}, both of which have been applied to bounds and synthesis problems involving pattern specification (\eg{}, directivity, pattern orthogonality) \cite{pozar1984antenna,LiuGarbaczPozar_AntennaSynthesisAndOptimizationUsingGeneralizedCharacteristicModes,Ethier_2009a}. Other modal problems may be constructed using various combinations of the total energy (rather than the reactance), Ohmic (thermal) losses, far-field, and radiation operators \cite{SchabBernhard_RadiationAndEnergyStorageCurrentModesOnConductingStructures, GustafssonTayliEhrenborgEtAl_AntennaCurrentOptimizationUsingMatlabAndCVX}, each with advantages and disadvantages leading themselves to specific problems, particularly in the study of bounds on antenna performance, \eg{},~\cite{jelinek2016optimal,EhrenborgGustafsson_FundamentalBoundsOnMIMOAntennas,Li_2019a}.  Similar in name to CMs, characteristic basis functions serve as a convenient method for reducing the computational cost of solving large systems by partitioning and expanding currents in terms of subdomain basis functions obtained for particular excitation \cite{prakash2003characteristic, Haykir2020}. 

The diagonalization properties of CMs make them attractive for accelerating optimization problems containing objectives or constraints involving the impedance matrix, such as those encountered in determining physical bounds on extinction and scattering cross sections \cite{2020_Gustafsson_NJP}. 

\section{Conclusion}
\label{sec:conclu}

For most of its history, CMA has been closely tied to numerical methods, particularly EFIE-based MoM for perfectly conducting structures, and since the publication of seminal work in the 1970's, continued research on CM decomposition has led to many findings in areas related to its computational implementation.  In this review, we touch on several of these areas, ranging from improved numerical techniques to new theoretical formulations.  Just as the field has progressed enormously in the last 50 years, we expect research in CMs to continue in directions yet-unknown for many years to come.

Perhaps just as impactful as CMA itself is the increased awareness and interest in source-based theoretical methods it has provided the antennas and propagation community.  The use of modal currents and fields as a basis for exploring the feasibility of electromagnetic devices is by no means unique to characteristic modes, though modern computing power has greatly accelerated progress in its application to the development of fast solution techniques, the calculation of physical bounds, and the discovery of further intuitions for practical antenna design.  Nevertheless, with every ``new finding'' in the 21$^\T{st}$ century, it seems the community also gains an increased appreciation for early pioneers of the field, \eg{}, the authors of work on numerical methods \cite{waterman1965matrix,Harrington_FieldComputationByMoM}, scattering decompositions \cite{1948_Montgomery_Principles_of_Microwave_Circuits, Garbacz_TCMdissertation} and source-based optimization \cite{lapaz1943optimum, uzsoky1956theory, chu1948physical, lo1966optimization,harrington1965antenna}; whose insights set the stage for modern work in CMA.

\bibliographystyle{IEEEtran}
\bibliography{IEEEabrv, nonIEEEabrv, references_CMA_review_project_no_URL_abbrev, references_CMA_review_Extra}

\begin{thebibliography}{10}
\providecommand{\url}[1]{#1}
\csname url@samestyle\endcsname
\providecommand{\newblock}{\relax}
\providecommand{\bibinfo}[2]{#2}
\providecommand{\BIBentrySTDinterwordspacing}{\spaceskip=0pt\relax}
\providecommand{\BIBentryALTinterwordstretchfactor}{4}
\providecommand{\BIBentryALTinterwordspacing}{\spaceskip=\fontdimen2\font plus
\BIBentryALTinterwordstretchfactor\fontdimen3\font minus
  \fontdimen4\font\relax}
\providecommand{\BIBforeignlanguage}[2]{{%
\expandafter\ifx\csname l@#1\endcsname\relax
\typeout{** WARNING: IEEEtran.bst: No hyphenation pattern has been}%
\typeout{** loaded for the language `#1'. Using the pattern for}%
\typeout{** the default language instead.}%
\else
\language=\csname l@#1\endcsname
\fi
#2}}
\providecommand{\BIBdecl}{\relax}
\BIBdecl

\bibitem{Harrington_FieldComputationByMoM}
R.~F. Harrington, \emph{Field Computation by Moment Methods}.\hskip 1em plus
  0.5em minus 0.4em\relax Piscataway, New Jersey, United States: Wiley -- IEEE
  Press, 1993.

\bibitem{Harrington_1971b}
R.~Harrington and J.~Mautz, ``Computation of characteristic modes for
  conducting bodies,'' \emph{{IEEE} Trans. Antennas Propag.}, vol.~19, no.~5,
  pp. 629--639, Sep 1971.

\bibitem{Makarov_AntennaAndEMModelingWithMatlab}
S.~N. Makarov, \emph{Antenna and EM Modeling with Matlab}.\hskip 1em plus 0.5em
  minus 0.4em\relax Wiley, 2002.

\bibitem{jin2011theory}
J.-M. Jin, \emph{Theory and computation of electromagnetic fields}.\hskip 1em
  plus 0.5em minus 0.4em\relax John Wiley \& Sons, 2011.

\bibitem{atom}
\BIBentryALTinterwordspacing
{A}ntenna {T}oolbox for {MATLAB} ({AToM}){,} {C}zech {T}echnical
  {U}niversity~in {P}rague, 2017. [Online]. Available:
  \url{www.antennatoolbox.com}
\BIBentrySTDinterwordspacing

\bibitem{feko}
\BIBentryALTinterwordspacing
A.~{FEKO}, Altair, 2017. [Online]. Available:
  \url{https://www.altair.com/feko/}
\BIBentrySTDinterwordspacing

\bibitem{Garbacz_TCMdissertation}
R.~J. Garbacz, ``A generalized expansion for radiated and scattered fields,''
  Ph.D. dissertation, The Ohio State Univ., 1968.

\bibitem{1948_Montgomery_Principles_of_Microwave_Circuits}
C.~G. Montgomery, R.~H. Dicke, and E.~M. Purcell, \emph{Principles of Microwave
  Circuits}.\hskip 1em plus 0.5em minus 0.4em\relax New York, United States:
  McGraw-Hill, 1948.

\bibitem{Harrington_1971a}
R.~Harrington and J.~Mautz, ``Theory of characteristic modes for conducting
  bodies,'' \emph{{IEEE} Trans. Antennas Propag.}, vol.~19, no.~5, pp.
  622--628, Sep 1971.

\bibitem{Harrington_1972b}
R.~Harrington, J.~Mautz, and Y.~Chang, ``Characteristic modes for dielectric
  and magnetic bodies,'' \emph{{IEEE} Trans. Antennas Propag.}, vol.~20, no.~2,
  pp. 194--198, Mar 1972.

\bibitem{SarkarMokoleSalazarPalma_AnExposeOnInternalResonancesCM}
T.~K. Sarkar, E.~L. Mokole, and M.~Salazar-Palma, ``An expose on internal
  resonance, external resonance and characteristic modes,'' \emph{{IEEE} Trans.
  Antennas Propag.}, vol.~64, no.~11, pp. 4695--4702, Nov. 2016.

\bibitem{Bernabeu_Jimenez_2017a}
T.~Bernabeu-Jimenez, A.~Valero-Nogueira, F.~Vico-Bondia, and A.~A. Kishk, ``A
  comparison between natural resonances and characteristic mode resonances of
  an infinite circular cylinder,'' \emph{{IEEE} Trans. Antennas Propag.},
  vol.~65, no.~5, pp. 2759--2763, May 2017.

\bibitem{Huang_StudyontheRelationshipsbetweenEigenmodesNaturalModesandCharacteristicModes}
S.~Huang, J.~Pan, and Y.~Luo, ``Study on the relationships between eigenmodes,
  natural modes, and characteristic modes of perfectly electric conducting
  bodies,'' \emph{International Journal of Antennas and Propagation}, vol.
  2018, p.~13, 2018.

\bibitem{RaoWiltonGlisson_ElectromagneticScatteringBySurfacesOfArbitraryShape}
S.~M. Rao, D.~R. Wilton, and A.~W. Glisson, ``Electromagnetic scattering by
  surfaces of arbitrary shape,'' \emph{{IEEE} Trans. Antennas Propag.},
  vol.~30, no.~3, pp. 409--418, May 1982.

\bibitem{Chang_1977a}
Y.~Chang and R.~Harrington, ``A surface formulation for characteristic modes of
  material bodies,'' \emph{{IEEE} Trans. Antennas Propag.}, vol.~25, no.~6, pp.
  789--795, Nov 1977.

\bibitem{Yla_Oijala_2019a}
P.~Yla-Oijala and S.~Jarvenpaa, ``Combined source integral equation-based
  theory of characteristic modes for impenetrable bodies,'' \emph{{IEEE} Trans.
  Antennas Propag.}, vol.~67, no.~4, pp. 2825--2828, Apr 2019.

\bibitem{Yla_Oijala_2019b}
P.~Yla-Oijala, ``Generalized theory of characteristic modes,'' \emph{{IEEE}
  Trans. Antennas Propag.}, vol.~67, no.~6, pp. 3915--3923, Jun 2019.

\bibitem{morse1953methods}
P.~M. {Morse} and H.~{Feshbach}, \emph{{Methods of theoretical physics}}, 1953.

\bibitem{Stratton_ElectromagneticTheory}
J.~A. Stratton, \emph{Electromagnetic Theory}.\hskip 1em plus 0.5em minus
  0.4em\relax Wiley -- IEEE Press, 2007.

\bibitem{Amendola1997}
G.~Amendola, G.~Angiulli, and G.~D. Massa, ``Numerical and analytical
  characteristic modes for conducting elliptic cylinders,'' \emph{Microw Opt
  Technol Lett.}, vol.~16, no.~4, pp. 243--249, Nov. 1997.

\bibitem{Capek_2017b}
M.~Capek, V.~Losenicky, L.~Jelinek, and M.~Gustafsson, ``Validating the
  characteristic modes solvers,'' \emph{{IEEE} Trans. Antennas Propag.},
  vol.~65, no.~8, pp. 4134--4145, Aug 2017.

\bibitem{Vandenbosch_ReactiveEnergiesImpedanceAndQFactorOfRadiatingStructures}
G.~A.~E. Vandenbosch, ``Reactive energies, impedance, and {Q} factor of
  radiating structures,'' \emph{{IEEE} Trans. Antennas Propag.}, vol.~58,
  no.~4, pp. 1112--1127, Apr. 2010.

\bibitem{GustafssonTayliEhrenborgEtAl_AntennaCurrentOptimizationUsingMatlabAndCVX}
M.~Gustafsson, D.~Tayli, C.~Ehrenborg, M.~Cismasu, and S.~Norbedo, ``Antenna
  current optimization using {MATLAB} and {CVX},'' \emph{{FERMAT}}, vol.~15,
  no.~5, pp. 1--29, May--June 2016.

\bibitem{Newman_SmallAntennaLocationSynthesisUsingCharacteristicModes}
E.~Newman, ``Small antenna location synthesis using characteristic modes,''
  \emph{{IEEE} Trans. Antennas Propag.}, vol.~27, no.~4, pp. 530--531, July
  1979.

\bibitem{Harrington_TimeHarmonicElmagField}
R.~F. Harrington, \emph{Time-Harmonic Electromagnetic Fields}, 2nd~ed.\hskip
  1em plus 0.5em minus 0.4em\relax Wiley -- IEEE Press, 2001.

\bibitem{PetersonRayMittra_ComputationalMethodsForElectromagnetics}
A.~F. Peterson, S.~L. Ray, and R.~Mittra, \emph{Computational Methods for
  Electromagnetics}.\hskip 1em plus 0.5em minus 0.4em\relax Wiley -- IEEE
  Press, 1998.

\bibitem{2018_Schab_Wsto}
K.~Schab \emph{et~al.}, ``Energy stored by radiating systems,'' \emph{{IEEE}
  Access}, vol.~6, pp. 10\,553--10\,568, 2018.

\bibitem{golub2013matrix}
G.~H. Golub and C.~F. Van~Loan, \emph{Matrix computations}.\hskip 1em plus
  0.5em minus 0.4em\relax JHU press, 2013, vol.~3.

\bibitem{matlab}
\BIBentryALTinterwordspacing
(2021) {MATLAB}. {The MathWorks}. [Online]. Available: \url{www.mathworks.com}
\BIBentrySTDinterwordspacing

\bibitem{Antonino_Daviu_2006a}
E.~Antonino-Daviu, C.~A. Suarez-Fajardo, M.~Cabedo-Fabres, and
  M.~Ferrando-Bataller, ``Wideband antenna for mobile terminals based on the
  handset {PCB} resonance,'' \emph{Microw Opt Technol Lett.}, vol.~48, no.~7,
  pp. 1408--1411, 2006.

\bibitem{Adams_2011a}
J.~J. Adams and J.~T. Bernhard, ``A modal approach to tuning and bandwidth
  enhancement of an electrically small antenna,'' \emph{{IEEE} Trans. Antennas
  Propag.}, vol.~59, no.~4, pp. 1085--1092, Apr 2011.

\bibitem{Ghosal_2020a}
S.~Ghosal, A.~De, A.~P. Duffy, and A.~Chakrabarty, ``Selection of dominant
  characteristic modes,'' \emph{{IEEE} Trans. Electromagn. Compat.}, vol.~62,
  no.~2, pp. 451--460, Apr 2020.

\bibitem{Tayli_2018a}
D.~Tayli, M.~Capek, L.~Akrou, V.~Losenicky, L.~Jelinek, and M.~Gustafsson,
  ``Accurate and efficient evaluation of characteristic modes,'' \emph{{IEEE}
  Trans. Antennas Propag.}, vol.~66, no.~12, pp. 7066--7075, Dec 2018.

\bibitem{Losenicky_etal_MoMandThybrid}
V.~Losenicky, L.~Jelinek, M.~Capek, and M.~Gustafsson, ``Method of moments and
  {T}-matrix hybrid,'' \emph{arXiv preprint arXiv:2012.04303}, 2020.

\bibitem{Dai_2017a}
Q.~I. Dai, H.~U.~I. Gan, Q.~S. Liu, and W.~C. Chew, ``Characteristic mode and
  reduced order modeling at low frequencies,'' \emph{IEEE Trans. Compon.
  Packag. Manuf. Technol.}, vol.~7, no.~5, pp. 669--677, May 2017.

\bibitem{yee1973self}
A.~Yee and R.~Garbacz, ``Self-and mutual-admittances of wire antennas in terms
  of characteristic modes,'' \emph{{IEEE} Trans. Antennas Propag.}, vol.~21,
  no.~6, pp. 868--871, 1973.

\bibitem{adams2013broadband}
J.~J. Adams and J.~T. Bernhard, ``Broadband equivalent circuit models for
  antenna impedances and fields using characteristic modes,'' \emph{{IEEE}
  Trans. Antennas Propag.}, vol.~61, no.~8, pp. 3985--3994, 2013.

\bibitem{vonNeumannWigner_UberDasVerhaltenEigenvertwen}
J.~von Neumann and E.~P. Wigner, ``\"{U}ber das verhalten von eigenwerten bei
  adiabatischen prozessen,'' \emph{Physicalische Zeitschrift}, vol.~30, pp.
  467--470, Sept. 1929.

\bibitem{Knorr_1973_TCM_symmetry}
J.~B. Knorr, ``Consequences of symmetry in the computation of characteristic
  modes for conducting bodies,'' \emph{{IEEE} Trans. Antennas Propag.},
  vol.~21, no.~6, pp. 899--902, Nov. 1973.

\bibitem{landau2013quantum}
L.~D. Landau and E.~M. Lifshitz, \emph{Quantum mechanics: non-relativistic
  theory}.\hskip 1em plus 0.5em minus 0.4em\relax Elsevier, 2013, vol.~3.

\bibitem{Schab_2017a}
K.~R. Schab and J.~T. Bernhard, ``A group theory rule for predicting eigenvalue
  crossings in characteristic mode analyses,'' \emph{{IEEE} Antennas Wireless
  Propag. Lett.}, vol.~16, pp. 944--947, 2017.

\bibitem{Masek_2020a}
M.~Masek, M.~Capek, L.~Jelinek, and K.~Schab, ``Modal tracking based on group
  theory,'' \emph{{IEEE} Trans. Antennas Propag.}, vol.~68, no.~2, pp.
  927--937, Feb 2020.

\bibitem{Raines_2012a}
B.~D. Raines and R.~G. Rojas, ``Wideband characteristic mode tracking,''
  \emph{{IEEE} Trans. Antennas Propag.}, vol.~60, no.~7, pp. 3537--3541, Jul
  2012.

\bibitem{Safin_2016a}
E.~Safin and D.~Manteuffel, ``Advanced eigenvalue tracking of characteristic
  modes,'' \emph{{IEEE} Trans. Antennas Propag.}, vol.~64, no.~7, pp.
  2628--2636, Jul 2016.

\bibitem{Chen_2021a}
X.~J. Chen, Y.~M. Pan, and G.~D. Su, ``An advanced
  eigenvector-correlation-based tracking method for characteristic modes,''
  \emph{{IEEE} Trans. Antennas Propag.}, vol.~69, no.~5, pp. 2751--2758, May
  2021.

\bibitem{Miers_2015a}
Z.~Miers and B.~K. Lau, ``Wideband characteristic mode tracking utilizing
  far-field patterns,'' \emph{{IEEE} Antennas Wireless Propag. Lett.}, vol.~14,
  pp. 1658--1661, 2015.

\bibitem{LudickJakobusVogel_AtrackingAlgorithmForTheEigenvectorsCalculatedWithCM}
D.~J. Ludick, U.~Jakobus, and M.~Vogel, ``A tracking algorithm for the
  eigenvectors calculated with characteristic mode analysis,'' in
  \emph{Proceedings of the 8th European Conference on Antennas and Propagation
  (EUCAP)}, 2014, pp. 569--572.

\bibitem{KingAJ_PhDThesis}
A.~J. King, ``Characteristic mode theory for closely spaced dipole arrays,''
  Ph.D. dissertation, University of Illinois at Urbana-Champaign, 2015.

\bibitem{Schab_2016a}
K.~R. Schab, J.~M. Outwater, M.~W. Young, and J.~T. Bernhard, ``Eigenvalue
  crossing avoidance in characteristic modes,'' \emph{{IEEE} Trans. Antennas
  Propag.}, vol.~64, no.~7, pp. 2617--2627, Jul 2016.

\bibitem{Lin_2017a}
J.-F. Lin and Q.-X. Chu, ``Extending bandwidth of antennas with coupling theory
  for characteristic modes,'' \emph{{IEEE} Access}, vol.~5, pp.
  22\,262--22\,271, 2017.

\bibitem{Ghosal_2020b}
S.~Ghosal, R.~Sinha, A.~De, A.~Chakrabarty, and H.~Son, ``Theory of coupled
  characteristic modes,'' \emph{{IEEE} Trans. Antennas Propag.}, vol.~68,
  no.~6, pp. 4677--4687, Jun 2020.

\bibitem{Vescovo_1989a}
R.~Vescovo, ``Characteristic modes for bodies endowed with mutually orthogonal
  symmetry planes,'' \emph{Microw Opt Technol Lett.}, vol.~2, no.~11, pp.
  390--393, Nov 1989.

\bibitem{Vescovo_1990a}
R.~Vescovo and T.~Corzani, ``Characteristic modes for nonconducting bodies
  having mutually orthogonal symmetry planes,'' \emph{Microw Opt Technol
  Lett.}, vol.~3, no.~4, pp. 124--127, Apr 1990.

\bibitem{Peitzmeier_2019a}
N.~Peitzmeier and D.~Manteuffel, ``Upper bounds and design guidelines for
  realizing uncorrelated ports on multimode antennas based on symmetry analysis
  of characteristic modes,'' \emph{{IEEE} Trans. Antennas Propag.}, vol.~67,
  no.~6, pp. 3902--3914, Jun 2019.

\bibitem{haus1991coupled}
H.~A. Haus and W.~Huang, ``Coupled-mode theory,'' \emph{Proc. {IEEE}}, vol.~79,
  no.~10, pp. 1505--1518, 1991.

\bibitem{Dai_2016a}
Q.~I. Dai, J.~Wu, H.~Gan, Q.~S. Liu, W.~C. Chew, and W.~E.~I. Sha,
  ``Large-scale characteristic mode analysis with fast multipole algorithms,''
  \emph{{IEEE} Trans. Antennas Propag.}, vol.~64, no.~7, pp. 2608--2616, Jul
  2016.

\bibitem{Wu_2018a}
B.-Y. Wu and X.-Q. Sheng, ``On the formulation of characteristic mode theory
  with fast multipole algorithms,'' \emph{{IEEE} Trans. Antennas Propag.},
  vol.~66, no.~11, pp. 6441--6445, Nov 2018.

\bibitem{Zhang_2020b}
Q.~Zhang, B.-Y. Wu, Y.~Gao, and X.-Q. Sheng, ``Multilevel fast multipole
  algorithm enhanced characteristic mode analysis for half-space platform,''
  \emph{{IEEE} Trans. Antennas Propag.}, vol.~68, no.~11, pp. 7711--7716, Nov
  2020.

\bibitem{Kabalan_2001a}
K.~Y. Kabalan, A.~El-Hajj, and A.~Rayes, ``A three-dimensional characteristic
  mode solution of two perforated parallel planes separating different
  dielectric media,'' \emph{Radio Science}, vol.~36, no.~2, pp. 183--193, Mar
  2001.

\bibitem{Liang_2018b}
P.~Liang and Q.~Wu, ``Duality principle of characteristic modes for the
  analysis and design of aperture antennas,'' \emph{{IEEE} Trans. Antennas
  Propag.}, vol.~66, no.~6, pp. 2807--2817, Jun 2018.

\bibitem{chalas2012efficient}
J.~Chalas and K.~Sertel, ``Efficient computation of in-situ antenna performance
  using platform characteristic modes,'' in \emph{Proceedings of the 2012 IEEE
  International Symposium on Antennas and Propagation}, 2012, pp. 1--2.

\bibitem{Ethier_2012a}
J.~Ethier and D.~McNamara, ``Sub-structure characteristic mode concept for
  antenna shape synthesis,'' \emph{Electronics Letters}, vol.~48, no.~9, p.
  471, 2012.

\bibitem{Ethier_2014a}
J.~L.~T. Ethier and D.~A. McNamara, ``Antenna shape synthesis without prior
  specification of the feedpoint locations,'' \emph{{IEEE} Trans. Antennas
  Propag.}, vol.~62, no.~10, pp. 4919--4934, Oct 2014.

\bibitem{Alroughani_2015a}
H.~Alroughani, J.~Ethier, and D.~A. McNamara, ``Orthogonality properties of
  sub-structure characteristic modes,'' \emph{Microw Opt Technol Lett.},
  vol.~58, no.~2, pp. 481--486, Dec 2015.

\bibitem{Alakhras_2020a}
A.~Alakhras and D.~A. McNamara, ``Sub-structure characteristic mode computation
  utilising field-based {MM}/{GTD} hybrid methods,'' \emph{Journal of
  Electromagnetic Waves and Applications}, vol.~34, no.~13, pp. 1812--1821, Jul
  2020.

\bibitem{Angiulli_2000a}
G.~Angiulli, G.~Amendola, and G.~Di~Massa, ``Application of characteristic
  modes to the analysis of scattering from microstrip antennas,'' \emph{Journal
  of Electromagnetic Waves and Applications}, vol.~14, no.~8, pp. 1063--1081,
  Jan 2000.

\bibitem{Borchardt_2019a}
J.~J. Borchardt and T.~C. Lapointe, ``U-slot patch antenna principle and design
  methodology using characteristic mode analysis and coupled mode theory,''
  \emph{{IEEE} Access}, vol.~7, pp. 109\,375--109\,385, 2019.

\bibitem{HarringtonMautz_PatternSynthesisForLoadedNportScatterers}
R.~F. Harrington and J.~R. Mautz, ``Pattern synthesis for loaded {N}-port
  scatterers,'' \emph{{IEEE} Trans. Antennas Propag.}, vol.~22, no.~2, pp.
  184--190, 1974.

\bibitem{Obeidat_2010b}
K.~A. Obeidat, B.~D. Raines, R.~G. Rojas, and B.~T. Strojny, ``Design of
  frequency reconfigurable antennas using the theory of network characteristic
  modes,'' \emph{{IEEE} Trans. Antennas Propag.}, vol.~58, no.~10, pp.
  3106--3113, Oct 2010.

\bibitem{Jaafar_2017a}
H.~Jaafar, S.~Collardey, and A.~Sharaiha, ``Optimized manipulation of the
  network characteristic modes for wideband small antenna matching,''
  \emph{{IEEE} Trans. Antennas Propag.}, vol.~65, no.~11, pp. 5757--5767, Nov
  2017.

\bibitem{Tzanidis_2012a}
I.~Tzanidis, K.~Sertel, and J.~L. Volakis, ``Characteristic excitation taper
  for ultrawideband tightly coupled antenna arrays,'' \emph{{IEEE} Trans.
  Antennas Propag.}, vol.~60, no.~4, pp. 1777--1784, Apr 2012.

\bibitem{Lonskyetal_CMofDipoleArrays}
T.~Lonsky, P.~Hazdra, and J.~Kracek, ``Characteristic modes of dipole arrays,''
  \emph{{IEEE} Antennas Wireless Propag. Lett.}, vol.~17, pp. 998--1001, June
  2018.

\bibitem{Guan_2019a}
L.~Guan, Z.~He, D.~Ding, and R.~Chen, ``Efficient characteristic mode analysis
  for radiation problems of antenna arrays,'' \emph{{IEEE} Trans. Antennas
  Propag.}, vol.~67, no.~1, pp. 199--206, Jan 2019.

\bibitem{Cheng_2019a}
G.~S. Cheng and C.-F. Wang, ``A novel periodic characteristic mode analysis
  method for large-scale finite arrays,'' \emph{{IEEE} Trans. Antennas
  Propag.}, vol.~67, no.~12, pp. 7637--7642, Dec 2019.

\bibitem{maalik2016characteristic}
A.~Maalik, R.~G. Rojas, and R.~J. Burkholder, ``Characteristic modal
  decomposition of reflection phase of a microstrip-patch reflectarray
  unit-cell,'' in \emph{2016 IEEE International Symposium on Antennas and
  Propagation (APSURSI)}, 2016, pp. 425--426.

\bibitem{Haykir_2019a}
Y.~Haykir and O.~A. Civi, ``Characteristic mode analysis of unit cells of
  metal-only infinite arrays,'' \emph{Advanced Electromagnetics}, vol.~8,
  no.~2, pp. 134--142, Sep 2019.

\bibitem{InagakiGarbacz_EigenfunctionsOfCompositeHermitianOperatorsWithApplicationToDiscreteAndContunuousRadiatingSystems}
N.~Inagaki and R.~J. Garbacz, ``Eigenfunctions of composite hermitian operators
  with application to discrete and continuous radiating systems,'' \emph{{IEEE}
  Trans. Antennas Propag.}, vol.~30, no.~4, pp. 571--575, July 1982.

\bibitem{Liu_SomeRelationshipsBetweenCharacteristicModesAndInagakiModesForUseInScatteringAndRadiationProblems}
D.~Liu, ``Some relationships between characteristic modes and {Inagaki} modes
  for use in scattering and radiation problems,'' Ph.D. dissertation, The Ohio
  State Univ., 1986.

\bibitem{pozar1984antenna}
D.~Pozar, ``Antenna synthesis and optimization using weighted {Inagaki}
  modes,'' \emph{{IEEE} Trans. Antennas Propag.}, vol.~32, no.~2, pp. 159--165,
  1984.

\bibitem{LiuGarbaczPozar_AntennaSynthesisAndOptimizationUsingGeneralizedCharacteristicModes}
D.~Liu, R.~J. Garbacz, and D.~M. Pozar, ``Antenna synthesis and optimization
  using generalized characteristic modes,'' \emph{{IEEE} Trans. Antennas
  Propag.}, vol.~38, no.~6, pp. 862--868, June 1990.

\bibitem{Ethier_2009a}
J.~Ethier and D.~McNamara, ``The use of generalized characteristic modes in the
  design of {MIMO} antennas,'' \emph{IEEE Transactions on Magnetics}, vol.~45,
  no.~3, pp. 1124--1127, Mar 2009.

\bibitem{SchabBernhard_RadiationAndEnergyStorageCurrentModesOnConductingStructures}
K.~R. Schab and J.~T. Bernhard, ``Radiation and energy storage current modes on
  conducting structures,'' \emph{{IEEE} Trans. Antennas Propag.}, vol.~63,
  no.~12, pp. 5601--5611, Dec. 2015.

\bibitem{jelinek2016optimal}
L.~Jelinek and M.~Capek, ``Optimal currents on arbitrarily shaped surfaces,''
  \emph{{IEEE} Trans. Antennas Propag.}, vol.~65, no.~1, pp. 329--341, 2016.

\bibitem{EhrenborgGustafsson_FundamentalBoundsOnMIMOAntennas}
C.~Ehrenborg and M.~Gustafsson, ``Fundamental bounds on {MIMO} antennas,''
  \emph{{IEEE} Antennas Wireless Propag. Lett.}, vol.~17, no.~1, pp. 21--24,
  jan 2018.

\bibitem{Li_2019a}
R.~Li, D.~A. McNamara, and G.~Wei, ``Evaluation of the available directivity of
  a radiating structure in terms of its characteristic mode content,''
  \emph{{IEEE} Trans. Antennas Propag.}, vol.~67, no.~10, pp. 6686--6691, Oct
  2019.

\bibitem{prakash2003characteristic}
V.~V.~S. Prakash and R.~Mittra, ``Characteristic basis function method: A new
  technique for efficient solution of method of moments matrix equations,''
  \emph{Microw Opt Technol Lett.}, vol.~36, no.~2, pp. 95--100, 2003.

\bibitem{Haykir2020}
Y.~Haykir and O.~A. Civi, ``Use of characteristic modes in the {CBFM} for the
  analysis of large arrays,'' in \emph{{EuCAP}}, Mar. 2020.

\bibitem{2020_Gustafsson_NJP}
M.~Gustafsson, K.~Schab, L.~Jelinek, and M.~Capek, ``Upper bounds on absorption
  and scattering,'' \emph{New Journal of Physics}, vol.~22, p. 073013, 2020.

\bibitem{waterman1965matrix}
P.~Waterman, ``Matrix formulation of electromagnetic scattering,'' \emph{Proc.
  {IEEE}}, vol.~53, no.~8, pp. 805--812, 1965.

\bibitem{lapaz1943optimum}
L.~La~Paz and G.~Miller, ``Optimum current distributions on vertical
  antennas,'' \emph{Proceedings of the IRE}, vol.~31, no.~5, pp. 214--232,
  1943.

\bibitem{uzsoky1956theory}
M.~Uzsoky and L.~Solymar, ``Theory of super-directive linear arrays,''
  \emph{Acta Physica Academiae Scientiarum Hungaricae}, vol.~6, no.~2, pp.
  185--205, 1956.

\bibitem{chu1948physical}
L.~J. Chu, ``Physical limitations of omni-directional antennas,'' \emph{Journal
  of Applied Physics}, vol.~19, no.~12, pp. 1163--1175, 1948.

\bibitem{lo1966optimization}
Y.~Lo, S.~Lee, and Q.~Lee, ``Optimization of directivity and signal-to-noise
  ratio of an arbitrary antenna array,'' \emph{Proc. {IEEE}}, vol.~54, no.~8,
  pp. 1033--1045, 1966.

\bibitem{harrington1965antenna}
R.~Harrington, ``Antenna excitation for maximum gain,'' \emph{{IEEE} Trans.
  Antennas Propag.}, vol.~13, no.~6, pp. 896--903, 1965.

\end{thebibliography}

\end{document}